
  \magnification=\magstep1
\settabs 18 \columns
\hsize=16truecm

\def\ss{\smallskip\smallskip}
\def\b{\bigskip}
\def\bb{\bigskip\bigskip}
\def\bbb{\bb\b}

\def\no{\noindent}
\def\r{\rightline}
\def\ce{\centerline}
\def\ve{\vfill\eject}

\def\r{\rightline}
  \font\got=eufm8 scaled\magstep1
  \def\g{{\got g}}

\def\harr#1#2{\smash{\mathop{\hbox to .25 in{\rightarrowfill}}
  \limits^{\scriptstyle#1}_{\scriptstyle#2}}}

\def\today{\ifcase\month\or January\or February\or March\or April\or
May\or June\or July\or
August\or September\or October\or November\or  December\fi
\space\number\day, \number\year }

\r \today
\bb\bb\bb

\def\Rit{\hbox{\it I\hskip -2pt R}}

\def\Cit{\hbox{\it l\hskip -5.5pt C\/}}
\def\Rrm{\hbox{\rm I\hskip -2pt R}}

\def\Crm{\hskip0.5mm \hbox{\rm l\hskip -4.5pt C\/}}

\def\sqr#1#2{{\vcenter{\vbox{\hrule height.#2pt
\hbox{\vrule width.#2pt height#2pt \kern#2pt
\vrule width.#2pt}
\hrule height.#2pt}}}}

  \def\1/2{{\scriptstyle{1\over 2}}}
  \def\a/2{{\scriptstyle{3\over 2}}}
  \def\5/2{{\scriptstyle{5\over 2}}}
  \def\7/2{{\scriptstyle{7\over 2}}}
  \def\3/4{{\scriptstyle{3\over 4}}}

\def\picture #1 by #2 (#3){
   \vbox to #2{
     \hrule width #1 height 0pt depth 0pt
     \vfill
     \special{picture #3} 
     }
   }

\def\scaledpicture #1 by #2 (#3 scaled #4){{
   \dimen0=#1 \dimen1=#2
   \divide\dimen0 by 1000 \multiply\dimen0 by #4
   \divide\dimen1 by 1000 \multiply\dimen1 by #4
   \picture \dimen0 by \dimen1 (#3 scaled #4)}
   }

%
%
\def\bbb{{\cal B}}

\font\steptwo=cmb10 scaled\magstep2

\magnification=\magstep1

   {\ce {\steptwo q-ALGEBRAS and ARRANGEMENTS OF HYPERPLANES}}
\b
\ce{Christian Fr\o nsdal}
\b
\ce{\it Physics Department, University of California, Los Angeles CA
90095-1547 USA}
\vskip1in

\def\sqr#1#2{{\vcenter{\vbox{\hrule height.#2pt
\hbox{\vrule width.#2pt height#2pt \kern#2pt
\vrule width.#2pt}
\hrule height.#2pt}}}}

\def \r{\rightarrow}

\no{\it ABSTRACT.} Varchenko's approach to quantum
groups, from the theory of arrangements of
hyperplanes, can be usefully applied to q-algebras in
general, of which
quantum groups and quantum
(super) Kac-Moody algebras are special cases. New
results are obtained on the classification of
q-algebras, and of the  Serre ideals of generalized quantum
(super) Kac-Moody algebras.
\b\b\b
\no {\steptwo 1. INTRODUCTION}
\b
\ce{\bf 1.1. Quantum Groups.}

Drinfel'd, in his address to the International Congress of
Mathematicians in Berkeley [D], defined what he proposed to
call quantum universal enveloping algebras, a class of
deformations of the enveloping algebras of Kac-Moody algebras. In
fact, the much wider family of generalized Kac-Moody algebras can
be similarly quantized, and there arises the new problem of
classifying these objects.  The generalized Kac-Moody algebras
themselves have  resisted classification till now, but because they
have a singular position within the deformed family (as is always
the case with essential deformations), there is some room for
hoping that an approach from general position may be
effective.

This section begins with a brief review of the structures
defined by Drinfel'd, setting the stage for introducing the
generalized  quantum groups and for a statement of the problem
addressed in this paper - in Subsection 1.1.4.
\ve

\no{\bf 1.1.1. Drinfel'd's Quantum Groups.}

Let \g ~ be a Kac-Moody algebra in the sense of Kac [K], defined in
terms of a `generalized Cartan matrix'. A square, complex
matrix  $A$ is so called if
$$
\eqalign{&
A_{ii} = 2,~~ i = 1,...,N,\cr &
A_{ij} {\rm ~is ~a ~non-positive ~integer~ for~} i\neq j,\cr &
A_{ij} = 0 {\rm ~implies~} A_{ji} = 0.\cr}
$$
Let $h$ be a complex parameter. For any generalized
Cartan matrix
$A$, the associated Drinfel'd Quantum Group (quantized Kac-Moody
algebra) is the
$\Crm[[h]]$-algebra generated by elements
$\{H_a,e_i,f_i\}^{a = 1,...,M}_{i = 1,...,N}$ with relations
$$\eqalign{&
[H_a,H_b] = 0, ~~a,b = 1,...,N,\cr &
[H_a,e_i] = H_a(i)e_i,~~[H_a,f_i] = -H_a(i)f_i,\cr &
[e_i,f_j] =  \delta_{ij}{2\over h}\sinh ({h\over 2} H^{\vee}
_i).\cr &
\cr}\eqno(1.1)
$$
\vskip-.4cm
\no Here $H(1),...,H(N)$ are the roots, $H_1^\vee,...,H_N^\vee$ are
the co-roots and $A_{ij} = H_i^\vee(H_j)$. Furthermore, for each
pair
$(i,j),	 i\neq j$, the quantum Serre relations
$$\eqalign{&
\sum_{k=0}^n (-1)^k\biggl({n\atop k}\biggr)_q
q^{-k(n-k)/2}(e_i)^ke_j(e_i)^{n-k}   = 0,\cr &
\sum_{k=0}^n (-1)^k\biggl({n\atop k}\biggr)_q
q^{-k(n-k)/2}(f_i)^kf_j(f_i)^{n-k}   = 0,\cr & q = {\rm e}^h,~~ n
= 1- A_{ij}.
\cr}\eqno(1.2)
$$
That $A$ is a generalized Cartan matrix implies that $n$ is a
positive integer. In [K], the author finds it convenient to
begin without this condition, taking $A$ to be an arbitrary matrix,
although ``a deep theory can be developed only for the Lie algebra
\g~associated to a generalized Cartan matrix ...". Our aim is to
challenge that remark. The question is what replaces the Serre
relations in the more general case.
\b

\no{\bf 1.1.2. Quantum Supergroups.}

The algebras defined in Subsection 1.1.1 are deformations of
Kac-Moody algebras. Super-Kac-Moody algebras can be
deformed in a similar manner, but the Serre relations are
more complicated and differ greatly from case to case. [Y]
The difference between Kac-Moody algebras and super-Kac-Moody
algebras has often been emphasized: they are different
types of tensor categories. But both categories merge upon
deformation; that is one of the attractive features of
quantization.
\b

\no{\bf 1.1.3. Serre relations.}

Let ${\cal A}$ be the algebra generated by
$\{H_a,e_i,f_i\}^{a = 1,...,M}_{i = 1,...,N}$ with relations
(1.1). Let ${\cal A}_+$ be the algebra generated by $\{e_i,H_a\}$
and the relations
$[H_a,H_b] = 0,~ [H_a,e_i] = H_a(i)e_i$ and let ${\cal A}_-$ be
the algebra generated by $\{f_i,H_a\}$ and the relations
$[H_a,H_b] = 0,~ [H_a,f_i] = -H_a(i)f_i$. Finally, let ${\cal
B}_+$ be the
\Crm-algebra freely generated by the $e_i$ and let ${\cal B}_-$
be the \Crm-algebra freely generated by the $f_i$.

Define mappings
$$\eqalign{ &
f_i^\#: {\cal B}_+ \rightarrow {\cal A}_+, ~e_j \mapsto
[e_j,f_i],\cr &
e_i^\#: {\cal B}_- \rightarrow {\cal A}_-,~f_j \mapsto
[f_j,e_i].\cr
\cr}
$$
An element of ${\cal B}_+~({\cal B}_-)$ is said to be
\underbar{invariant} if it is annihilated by all the mappings
$f_i^\#$  (all the mappings $e_i^\#$). Let ${\cal I}_+ \subset
{\cal B}_+$ and $  {\cal I}_-
\subset {\cal B}_-$ be the two-sided ideals generated by the
invariants.
\b

\no{\bf Definition.} {\it The Serre ideal of ${\cal A}$ is the
direct sum ${\cal I}({\cal A}) = {\cal I}_+ \oplus {\cal I}_-$.}
\b
\no{\bf Theorem.} {\it The Serre ideal of ${\cal A}$ is
generated by the Serre relations.}
\b

This allows to define the Drinfel'd quantum group as the algebra
$$
{\cal A}' = {\cal A}/{\cal I}({\cal A}),\eqno(1.3)
$$
and this formulation allows to relax the condition that $A$
be a generalized Cartan matrix, and to define a generalized
quantum group.
\b

\no{\bf 1.1.4. Generalized quantum groups.}
\b
\no{\bf Definition.} {\it Let ${\cal M} ,{\cal N}$ be two
countable sets, and
$\phi,\psi$ two maps,
$$\eqalign{&
\phi: {\cal M}\times {\cal M} \rightarrow \Crm,~~a,b \mapsto
\phi^{ab},\cr &
\psi: {\cal M}\times {\cal N} \rightarrow \Crm, ~~a,i \mapsto
H_a(i).
\cr}
$$
Let
$$
\phi(i,\cdot) =
\sum_{a,b\in {\cal M}}\phi^{ab}H_a(i)H_b,~\phi(\cdot,i) =
\sum_{a,b\in {\cal M}}\phi^{ab}H_aH_b(i),
$$
and suppose that ${\rm e}^{\phi(i,\cdot) + \phi(\cdot,i)}
\neq 1,
i\in {\cal N}$.

Let ${\cal A} = {\cal A}(\phi,\psi)$ be the universal, associative,
unital
\Crm-algebra with generators}

\no{\it $\{e_i,f_i\}_{i\in
{\cal N}}$ and
$\{H_a\}_{a\in {\cal M}}$ and relations
$$\eqalign{&
[H_a,H_b] = 0,~~ a,b\in {\cal M},\cr &
[H_a,e_i] = H_a(i)e_i,~~[H_a,f_i] =
-H_a(i)f_\beta,\cr &
[e_i,f_j] = \delta_{ij}\bigl({\rm
e}^{\phi(i,\cdot)} - {\rm e}^{-\phi(\cdot,i)}\bigr).\cr}
$$
Then the generalized quantum group ${\cal A}'={\cal A}'(\phi,\psi)$
is the quotient ${\cal A}' = {\cal A}/{\cal  I}({\cal A})$, where
${\cal I}({\cal A})$ is the Serre ideal of
${\cal A}$, defined as in Subsection 1.1.4.}
\b
If the form $\phi$ is symmetric and $A_{ij} = \phi(i,j)$ is
a generalized Cartan matrix, then ${\cal A}'$ is a quantum
group in the sense that this term is used in most of the
litterature. The partial generalization
that consists of relaxing the symmetry requirement was studied
by Reshetikhin [Re]. The program of this paper is the
classification  of the larger family of generalized quantum groups
(no restrictions on the matrix $A$),  in terms of their Serre
ideals.
\b
\no{\bf 1.1.5. Generalized Drinfel'd-Jimbo algebras.}

Technical difficulties that arise from the appearance of infinite
series in $H_1,...,H_M$ within the relations can be
avoided.  If $M = N$ and $H_i^\vee = 2H_i$ replace the
Cartan generators $H_1,...,H_N$ by
$$
K_i = {\rm e}^{H_i}, ~K^i = {\rm e}^{-H_i},~~ i = 1,...,N~.
$$
The relations are now
$$\eqalign{&
K_iK^i = K^iK_i =1,\cr &
K_ie_jK^i = {\rm e}^{H_i(j)}e_j, ~~K_if_jK^i = {\rm
e}^{-H_i(j)}f_j,\cr & [e_i,f_j] = \delta_{ij}{2\over
h}(K_i-K^i),\cr}\eqno
$$
and the Serre relations. A slight disadvantage is that the
classical limit is no longer the underlying Kac-Moody algebra \g~;
the difference arises from the fact that $K_i \rightarrow \pm 1$.
However, since statements that are true for Drinfel'd's quantized
enveloping algebras usually imply analogous results for the
Drinfel'd-Jimbo algebras we shall not be greatly concerned with
the distinction.

In the general case set
  $$
K_i = {\rm e}^{\phi(i,.)},~~ K^i = {\rm e}^{-\phi(.,i)}.
$$
The relations $K_iK^i = K^iK_i = 1$ are omitted while the rest of
the relations remain as written.
\b
\no{\bf 1.1.6. Generalized quantum supergroups.}

An interesting aspect of super Lie algeras is the existence of two
kinds of odd roots. In the case of a generalized quantum
supergroups there are parameters $q_{ii}$ that are fixed and equal
to -1; then $e_i$ is a `null root' and one of the relations is
$e_i^2 = 0$. The other odd roots are characterized by the fact
that $q_{jj}\rightarrow -1$ in the classical limit.

\b\b
\ce{\bf 1.2. Hopf structure of generalized quantum groups.}

In Drinfel'd's terminology a quantum group is a coboundary Hopf
algebra. This Hopf structure plays a relatively minor role in this
paper since the methods used are essentially algebraic. However,
the differential operators $\partial_i$  (see
Subsection 1.2.4.)  first appeared in an investigation of Hopf
structures; this justifies a short review. It is possible that
a more direct use of the Hopf structure may lead to simpler proofs
and, in the hands of an expert, to further results.
\ve
\no{\bf 1.2.1. Hopf Structure.}

Fix the sets ${\cal M},{\cal N}$ and the maps $\phi,\psi$ and let
${\cal A},  {\cal I}$ and ${\cal A}'$ be defined as in 1.1.4.
\b
\no{\bf Proposition.}[F1]
{\it There exists a unique homomorphism $\Delta: {\cal A}
\rightarrow {\cal A}\otimes {\cal A}$, such that
$$\eqalign{ &
\Delta(H_a) = H_a\otimes 1 + 1\otimes H_a, ~~A\in {\cal M},\cr &
\Delta(e_i) = 1\otimes e_i + e_i \otimes
{\rm}^{\phi(i,\cdot)},\cr &
\Delta(f_i) = {\rm e}^{-\phi(\cdot,i)}\otimes f_i +
f_i \otimes 1,~~ i \in {\cal N}.\cr}
$$
   The homomorphism $\Delta$ induces a unique homomorphism $
{\cal A}' \rightarrow A'\otimes {\cal A}'$, also denoted $\Delta$.
The algebra ${\cal A}$ becomes a Hopf algebra when endowed with
the counit ${\cal E}$ and the antipode $S$. The former is
the unique homomorphism ${\cal A}  \rightarrow {\cal A} $
that vanishes on all the generators. The antipode is the unique
anti-homomorphism ${\cal A} \rightarrow {\cal A} $ such that
$$\eqalign{&
S(H_a) = -H_a,~~ a \in {\cal M},\cr &
S(e_i) = - e_i {\rm e}^{-\phi(i,\cdot)},~~
S(f_i) = -{\rm e}^{\phi(\cdot,i)}f_i,~~i \in
{\cal N}.
\cr}
$$
The counit ${\cal E}$ and the antipode $S$ induce analogous
structures on ${\cal A}'$.}
\b
\no{\bf 1.2.2.  Coboundary property.}

Let $\Delta'$ denote the opposite coproduct: in Sweedler's
notation, if
$\Delta(x) =
\sum x^1\otimes x^2$ then $\Delta'(x) = \sum x^2\otimes x^1$.

We restrict our attention temporarily to the special case of
Drinfel'd's quantum groups. Then there exists an
element $R \in {\cal A}'\otimes {\cal A}'$ that interpolates
between $\Delta$ and $\Delta'$:
$$
\Delta(x)R - R\Delta'(x) = 0,~~ \forall x \in {\cal A}'.
$$
This element is known as the Universal Yang-Baxter
Matrix; it satisfies the Yang-Baxter relation
$$
R_{12}R_{13}R_{23} = R_{23}R_{13}R_{12},
$$
and it has been calculated explicitly by Reshetikhin [Re] and
others.

Let $\sigma$ be the operator in ${\cal A}'\otimes {\cal A}'$ that
interchanges the two spaces, and let $P := \sigma \circ R$. Then
$$
Q(x) := \Delta(x)P-P\Delta(x) = 0.
$$
Let $ d:~ $Hom$({\cal A}'^{\otimes p},{\cal A}'^{\otimes
q})\rightarrow ~$Hom$({\cal A'}^{\otimes p+1},{\cal A'}^{\otimes
q })$ be the Hochschild differential of ${\cal A}'$. We have $P\in
~$Hom$ (\Crm,{\cal A}'^ {\otimes
2})$, and $d P = Q\in~ $Hom$({\cal A}',{\cal A}'^{\otimes
2})$. Let
$U$ be the bialgebra topologically dual to ${\cal A}'$.
By duality, $Q(x)$ is interpreted as an element of   Hom
$(U^{\otimes 2}, \Crm)$ and $Q(x) = 0$ determines the algebraic
structure of
$U$.   We have
  $$
dQ(x,y) = \Delta(x) Q(y) - Q(xy) +
Q(x)\Delta(y)\in {\rm Hom} (U^{\otimes 2},\Crm).
$$
If $Q(x)  = Q(y)  =0,~ x,y \in {\cal A}' $, then the property
$d Q = 0$ reduces to $Q(xy)  =
0$. Thus
$Q$ must be closed, while the existence of the Universal R-matrix
tells us that $Q$ is exact. Accordingly, the Drinfel'd quantum
groups are called coboundary Hopf algebras.

\b

\no{\bf 1.2.3. The R-matrix of a generalized quantum group.}

The Hopf algebras introduced in Subsection 1.1.4 are also of the
coboundary type, and the R-matrices have been calculated in [F1].
\b
\no{\bf Proposition.} [F1] {\it The algebra ${\cal A}'$ is a
coboundary Hopf algebra with a Universal R-matrix in the form of a
series
$$
R = {\rm e}^\phi\bigl(1 + f_i\otimes e_i + \sum_{n =
2}^\infty t_{\underline i}^{\underline j}~f_{\underline
i}\otimes e_{\underline j}\bigr),
$$
with
$$
\phi = \phi^{ab}H_a\otimes H_b,~\underline i = i_1,...,i_n,
~f_{\underline i} = f_{i_1}  ...
f_{i_n},
$$
and with complex coefficients $t_{\underline i}^{\underline
j}$. }
\b
\no{\bf Outline of proof.}
(a) Define elements $t_{\underline i}\in {\cal B}_+$ by
$$
t_{\underline i} = \sum_{\underline
j}t_{\underline i}^{\underline
j}e_{\underline j}.
$$
  By direct calculation one finds that the above series satisfies
the Yang-Baxter relation if and only if the following recursion
relations hold,
$$
[t_{\underline i},f_k] = {\rm
e}^{\phi(k,\cdot)}\delta_{i_1}^k
t_{\underline i} -
t_{\underline i}\delta_{i_n}^k
e^{-\phi(\cdot,k)},~~ k = 1,2,...,N.
$$

\no (b)
  Define operators $\partial_k$ on ${\cal B}_+$ by
$$
\partial_k e_i x = \delta_{ik} x +
e^{-\phi(k,i)}e_i \partial_k x,~~x\in {\cal
A}'.
$$
then the above recursion relation is equivalent to
$$
\partial_k t_{\underline i} = \delta_{i_1}^k
t_{\underline i}.
$$

\no (c) Define the matrix $S$ by
$$
S_{\underline i}^{\underline j} = \partial_{\underline
j} ~e_{\underline i},~~ \partial_{\underline j} =
\partial_{j_n}...\partial_{j_1},\eqno(1.4)
$$
for multi-indices of equal length $ n = 1,2,...~$, all other
matrix elements zero. The
   projection of this matrix on ${\cal A}'$ is invertible, and the
inverse is the projection on ${\cal A}'$ of the matrix $t$ with
matrix elements $t_{\underline i}^{\underline j}$.

\no(d) Finally it is easy to verify that this R-matrix satisfies
the  relation $\Delta(x)R = R\Delta'(x)$.

\b
The complete proof makes extensive use of the properties of the
algebra
${\cal B}_+$ endowed with the differential structure introduced
by the action of the operators $\partial_k$. Some of these
properties will be summarized below.
\b

\no{\bf 1.2.4.} {\bf Proposition.}

{\it The space ${\cal I}_+
\subset {\cal B}_+$, generated by the invariants in
${\cal B}_+$, coincides with the space generated by the
``constants''; namely,   the elements $x\in {\cal B}_+$ that
satisfy the relations
$
\partial_i x = 0, ~i = 1,...,N~.
$}

\b
The problem of determining the Serre ideal of ${\cal A}'$ is thus
reduced to the calculation of the space of constants in ${\cal
B}_+$.
\b

\ce{\bf 1.3. Classification of q-algebras.}

It is proposed to determine the Serre ideals of the algebras
${\cal A}' ={\cal
A}'(\phi,\psi)$ defined in Subsection 1.1.5. These algebras are
parameterized by the values of the maps $\phi$ and $\psi$, the
ideals by the parameters
$$
q_{ij} = {\rm e}^{-\phi(i,j)}\neq 0,~~i,j = 1,...,N~.
$$
By Proposition 1.2.4 the problem reduces to a study of
q-algebras, that we now define.
\b

\no{\bf 1.3.1. Definition; q-algebras.}

{\it On the freely generated
algebra $\bbb = \Cit[e_1,\cdots,e_N]$, introduce
differential operators $\partial_1,\cdots,\partial_N$
with the action defined by $\partial_ie_j =
\delta_{ij}$ and
$$
\partial_i(e_jx) = \delta_{ij} x + q_{ij} e_j\partial_i
x, \quad x \in \bbb.
$$
Let $\bbb_q$ be the same algebra $\bbb$ with this
differental structure.}
\b
\no{\bf 1.3.2. Definition; constants.}

{\it  A `constant' in $\bbb_q$  is
a polynomial $C\in \bbb$, having no term of order 0,
such that $\partial_iC = 0,~ i = 1,\cdots,N$. Let ${\cal
I}_q$ denote the ideal in $\bbb_q$ that is generated by
the constants.}
\b

\no{\bf Theorem.} ([F1])  {\it The ideal ${\cal I}_q$ of ${\cal
B}$, via the identification of ${\cal B}$ with the subalgebra
${\cal B}_+\subset {\cal A}$, is precisely the component
${\cal I}_+$ of the Serre
ideal of ${\cal A}$.}
\b
The interest focuses on the quotient,
$$
\bbb_q' := \bbb_q/{\cal I}_q.\eqno(1.5)
$$
These are the `q-algebras' of the title.
\b
\no{\bf 1.3.3. Remarks.}

(a) One can introduce a second set of
differential operators $\partial'_i$, acting on $\bbb$
from the right. Kharchenko has shown [Kh1] that
$(\bbb,\partial_i,\partial'_i)$ is a bicovariant
differential structure in the sense of Woronowicz [W].
(b) The operators $\partial_i$ and $\partial'_i$ have
been introduced by Kashiwara, in his work [Ka] on crystal
bases. Kashiwara showed that these operators can be constructed
inside the quantum group (with one parameter). Whether this
remains true for generalized quantum groups is not known,
nor is it directly relevant for the problematics of this
paper.
\b
\no{\bf 1.3.4. Essential parameters.}

   It is shown in [FG] that the essential parameters,
that determine the existence and the coefficients of constants, are
$$
\sigma_{ij} := q_{ij}q_{ji}, \quad i\neq j,\quad i,j =
1,\cdots N,
$$
and $q_{ii},~ i = 1,\cdots,N$. For generic values of
these parameters there are no constants and the Serre ideal
of ${\cal A}$ is empty.

\b
\no{\bf 1.3.5. Gradings.}

The algebra $\bbb$ has a natural
grading by the total polynomial degree, and this grading
is passed on to $\bbb_q$  and to
$\bbb_q'$. A finer grading is the map that takes $e_{i_1}
\cdots e_{i_n} $ to the unordered set
$\{i_1,\cdots,i_n\}$. Under this grading, the monomials
of $\bbb$ (or $\bbb_q$ or $\bbb_q'$) are partially
ordered by the relation of inclusion of sets,  which gives a
sense to
   the term `lower
degree'. The space of constants has a basis of polynomials
that are homogeneous in this finer grading; that is,
linear combinations of the permutations of a single
monomial.
\b
\no{\bf Definition.} {\it A homogeneous constant is called
`primitive' if it is not in the ideal generated by constants of
lower degree. The `space of primitive constants' is defined
via filtration.}
\b
\no The ideal ${\cal I}_q$ is generated by a set of
  primitive constants.

\b
The first general result was this.
\b
  \no{\bf 1.3.6. Theorem.}  [FG]) {\it Fix the degree
$G =\{1,\cdots,n\}$, and suppose that there are no
constants of lower degree. Then the space of constants
of   degree $G$ has dimension
$$
\eqalign{&
(n-2)!,\quad {\rm if }~\sigma_{1\cdots n} := \prod
q_{ij} = 1,\cr
&  \quad \quad \quad 0,\quad {\it otherwise}.\cr}
$$
The product runs over all pairs $i\neq j,~ i,j =
1,\cdots,n$.}
\b
\no An essentially equivalent result was obtained by
Kharchenko [Kh2].
\b
\no {\bf 1.3.7. Comparison with other work.}

  In addition to
references already quoted we mention the work of Rosso [R]. He
gives a nice direct presentation of the $q$-algebras ${\cal B}'_q$
in which the ideal vanishes identically. This is equivalent to a
result in [FG], where it was shown that the homomorphism from
${\cal B}_q$ to the algebra
${\cal B}_q^*$ of quantum differential operators on ${\cal B}_q$,
defined by
$e_i\mapsto\partial _i$, induces an isomorphism between ${\cal
B}'_q$ and ${\cal B}_q^*$. The Hopf structure is
prevalent in the work of  Flores de Chelia and Greene [FCG], who
have recently arrived independently at a result that is
equivalent to Theorem 1.3.6. In our work
the Hopf structure is represented by the matrix $S$ (Section
2).
This matrix is intimately related to the universal R-matrix
(Subsection 1.2.3); it appears in almost all work in this area,
notably in Varchenko [V] (who regards it as a form and calls it
$B$) and in the paper  [FCG] (where it is denoted $\Omega$).

\ve
\ce{\bf 1.4. Summary.}

General results for the case
of arbitrary degree $G$, but with the essential proviso
that there be no constants of lower degree, have been
reported [F2]. In this paper we return to the multilinear case $G =
\{1,\cdots,n\}$ (no repetitions).
In Section 2 we reduce the problem to a study of a
determinant, and set up a scheme for the classification of
q-algebras in terms of determinatal varieties. In Section 3 we
explain the results of Varchenko that will be used. In Sections 4
and 5 we limit our study to
  the case when there may be any number of
primitive constants of lower order, but all of total
degree 2.
These constants are generated by polynomials of the form
$$
e_ie_j - q_{ji}e_je_i,\quad i\neq j,\eqno(1.4)
$$
and these polynomials are constants if and only if
$\sigma_{ij} := q_{ij}q_{ji} = 1$.

Results for this special case are obtained
in Section 4 and  presented as Theorem 4.2. As I do
not know how or if the method of arrangements of
hyperplanes can be adapted to a more general
situation, I present in Section 5 an alternative and
completely algebraic proof of Theorem 4.2. Though it
owes much to the paper [V], it makes no use of
geometric concepts.

In Section 6 it is shown that this
new approach is applicable to a much more general
case, allowing for any number and any type of
constraints (and constants) of lower degree. The
result, Theorem 6.5, is a  solution for the
multilinear case, $G = \{1,\cdots ,n\}$, under the
stipulation that there be at least one pair $\{i,j\}$,
such that there is no constraint on $\sigma_{ij}$.
This last stipulation is important; unfortunately it is
violated by ordinary quantum (super) groups.

  \bb
\no{\steptwo 2. The matrix S and the form B.}
\b
\ce{\bf 2.1. The matrix S.}
We continue to use the multi index
notation,
$\underline i := i_1\cdots i_n, ~~\underline i' =
i_n\cdots i_1$ and
$$
\partial_{\underline i'} = \partial_{i_n}\cdots
\partial_{i_1},\quad e_{\underline j} = e_{j_1}\cdots
e_{j_n}.
$$
A matrix $S = (S_{\underline i\underline j})$ is defined
by
$$
   S_{\underline i}^{\underline j}  = \partial_{\underline j'}
e_{\underline i} |_0,
$$
where $x|_0$ is the term of total order 0 in the
polynomial $x \in \bbb$.
This matrix commutes with the grading,
$$
S =  \oplus_G S_G,\quad  (S_G)_{\underline i\underline
j} = \partial_{\underline i'} e_{\underline j},\quad
$$
where $\underline i, \underline j$ run over the
orderings of the unordered set $G$.

The matrix $S$ is singular if and only if there is a
constant in $\bbb_q$, and $S_G$ is singular if and only
if there is a constant (primitive or not) of degree $G$.
The existence of constants can thus be decided by
inspection of the determinants. For example, if~
$\sigma_{12} := q_{12}q_{21} = 1$, then there is a
constant of degree $G = \{1,2\}$, namely
$e_1e_2 - q_{21}e_2e_1$, and
$$
S_G = \pmatrix{1&q_{12}\cr q_{21} & 1},\quad \det S_G =
1-\sigma_{12} = 0.
$$
\ve

\ce{\bf 2.2.   The determinant.}

\no{\bf 2.2.1. Parameters in general position.}

The family $\{{\cal B}_q\}$ of
algebras is parameterized by
$q = \{q_{ij}\}_{i,j =
1,...,N}\in V :=\Crm^{N^2}$. There is an open subset $V_{gen}$
of $V$ such that for $q\in V_{gen}$ there are no
constants in ${\cal B}_q$, namely, the subspace defined by $\det S
\neq 0$. We shall say that parameters in this open set are in
general position. Until further notice suppose that the parameters
are in general position.

Let $\bbb_G$ be the
subspace of $\bbb_q$ that consists of all polymials of
degree G. From now on in this paper $G = \{1,\cdots ,n\}, ~
n$ fixed.   Set
$$
w_{n,k} = u_{n,k}v_k~,\eqno(2.1)
$$
where
$$
u_{n,k} = (n+1-k)!\eqno(2.2)
$$
  and
$$
\quad v_k =
(k-2)! \eqno(2.3)
$$
Then it is a result of Varchenko that
$$
\det S_G =
\prod_k\prod_{i_1,\cdots ,i_k}(1-\sigma_{i_1\cdots
i_k})^{w_{n,k}}.\eqno(2.4)
$$
The inner product is over all subsets of cardinality $k\geq 2$
of the set $\{1,\cdots ,n\}$.
  The total degree   in $q$'s  of $\det S_G$ is $\pmatrix{
n \cr 2 \cr} n! $~, and the formula implies the sum rule
$$
\sum_{k=2}^nk(k-1)w_{n,k}\pmatrix{n\cr k\cr} =
\pmatrix{n\cr 2}n! \eqno(2.5).
$$
Since all $\sigma_{ij}$ appear symmetrically, the total
degree in $\sigma_{12}$, say, is
$$
\sum_{k=2}^n w_{n,k}\pmatrix{n-2\cr k-2\cr} =
n!/2 .\eqno(2.6)
$$
The numbers (2.2) and (2.3) have the following interpretation.
Fix the integer $k \leq n$ and let $G_k = \{1,...,k\}$.  Let the
parameters approach a portion of the boundary of $V_{gen}$  where
$\sigma_{1...k} = 1$ but
$\sigma_{\underline i} \neq 1$ for all $\underline i \neq 1...k$
(as un-ordered sets). Then primitive constants appear in
${\cal B}_{G_k}$;  $v_k$ is the dimension   of the
space of (primitive) constants in ${\cal B}_{G_k}$ and
$u_{n,k}$ is the dimension of the ideal in ${\cal B}_G$ generated
by each constant in ${\cal B}_{G_k}$. A geometrical interpretation
will follow.

\b
\no{\bf  2.2.2. Example.}

  Let $G = \{1,2,3\}$ and suppose that
there are no constants of lower degree, then
$$
\det S_G =
(1-\sigma_{12})^2(1-\sigma_{23})^2(1-\sigma_{13})^2(1-\sigma_{123}).
$$
The surface on which $S_G$ is singular has four components, and
in particular $S_G$ is singular on the surface $\sigma_{123} = 1$.
On this surface the algebra ${\cal B}_q$ is characterized by the
existence of a primitive constant of degree $G = \{123\}$.
\b
\no{\bf  2.2.3. Example.}

Let $G = \{1,2,3,4\}$ and suppose that
there are no constants of lower degree, then
$$
S_G = \prod_{i<j}(1-\sigma_{ij})^6\prod_{i<j<k}(1-\sigma_{ijk})^2
(1-\sigma_{1234})^2.
$$
On the surface $\sigma_{1234} = 1$ there is a 2-dimensional
subspace of constants in ${\cal B}_G$.

\b
\ce{\bf 2.3. Cell decomposition of parameter space.}

The space of parameters is the space $V = \Crm^{N^2}$ in which the
$N^2$ parameters $q_{ij}$ take their values, with the natural
analytic structure defined by these parameters. This space is the
disjoint union of its $G$-cells ($G$ fixed),
  defined as follows.
\b

\no{\bf 2.3.1. Definition.}

{\it A $G$-cell in $V$ is a connected
subset of $V$ on which the rank of each matrix $S_{G'},~
G'\leq G,$ is   constant. A  regular function on a $G$-cell is the
restriction to the cell of a polynomial on $V$.} \break
  \vskip-3mm
   There is a space of constants associated to each  point   $q\in
V$,  and a regular field of constants on each $G$-cell.

\b
\no{\bf 2.3.2. Definition.}

{\it Two algebras ${\cal B}'_q$ and ${\cal
B}'_{q'}$ are of the same
$G$-type if $q$ and $q'$ belong to the same $G$-cell. They are of
the same multilinear type if they are of the same
$G$-type for every degree $G = \{i,j,...\}$ without repetition, and
of the same type if they are of the same $G$-type for every degree
$G$.}
\b
\no{\bf 2.3.3. Classification of the algebras by type.}

It is our
final aim to classify
$q$-algebras by type; in this paper we have the more modest goal
of a preliminary classification by multilinear type. The proposed
strategy is inductive. For $ G = \{1,2\}$ there are are two
cells:
$$
C_1: \sigma_{12} \neq 1,~~C_2 = dC_1: \sigma_{12} = 1,
$$
where $dC$ denotes the boundary of the cell $C$. Suppose the cells
have been determined for all multilinear degrees lower than $G =
\{1,2,...,n\}$. Fix the $G'$-type for each $G'<G$; this amounts
to fixing a certain set $Q$ of constraints of lower order.  Let
$V^Q$ be the closed subspace of $V$ defined by these constraints,
let ${\cal I}(Q)$ be the  ideal generated by the associated
constants, ${\cal B}_G(Q)=
{\cal B}_G/I(Q)\cap {\cal B}_G$ and $S_G(Q)$ the projection of
$S_G$ on ${\cal B}_G(Q)$. There is an open subset of
$V^Q$ on which
$\det S_G(Q)
\neq 0$ and ${\cal B}_G(Q)$ has no constant, and this
determines the G-type of
${\cal B}_q'$ for these parameter values. There remains the
boundary
$dV^Q$ of
$V^Q$, the hypersurface $\det S_G(Q) = 0$. The points of this
boundary are of two kinds. First, those characterized by the
appearance of one or more additional constant of lower
degree, each determined by a constraint that involves a proper
subset of
$\{\sigma_{ij}\}_{i,j = 1,...,n}$. This places the parameter in a
$V^{Q'}$ of lower dimension. By treating the spaces $V^Q$ in order
of non-increasing dimension we avoid having to take these points
into  account.
  The complement in $dV^Q$ of this first part of the
boundary, if any,  will be called the `primitive boundary
of $V^Q$;, it consists of points where
$\sigma_{1...n} = 1$.  The classification of types reduces to the
question of the existence of primitive boundaries.

In Example 2.2.2 above the set $Q$ is empty and $S_G(Q) = S_G$.
The expression for the determinant shows that there is a
primitive boundary characterized by the constraint
$ \sigma_{123} = 1$ (and $\sigma_{12},\sigma_{13},\sigma_{23} \neq
1$). The
$G$-cells $C_1,C_2$  are the subsets of $V^Q$ defined by
$\sigma_{12},\sigma_{13},\sigma_{23} \neq 1$ and
$$
C_1: ~ \sigma_{123}\neq 1, ~~  ~C_2:  ~ \sigma_{123} =
1.
$$
This situation is further illustrated by Example 2.2.3. Here too
there are two $G$-cells of interest, on which all $\sigma_{ij}\neq
1$, all $\sigma_{ijk} \neq 0$ and $\sigma_{1234}$ is either equal
to 0 or  different from zero.

\b

\no{\bf 2.3.3.  Example.}

Let $Q$ be the constraint $\sigma_{12} = 1$ associated with the
constant $e_1e_2 - q_{21} e_2e_1$ and
$G =
\{1,2,3\}$. Then
$V^Q$ is the surface in $V$  on which $\sigma_{12} = 1$, and
$$
\det S_G(Q) = (1-\sigma_{13})^2(1-\sigma_{23})^2\neq 0.
$$
There is no primitive boundary and only one classifying $G$-type in
this case. In the generic case  a new constant appears on the
surface
$\sigma_{123} = 1$, but in the present special case, when there is
a constant of lower order, this
surface is not singular for $S_G(Q)$.

\b\no{\bf  2.3.4. Example.}

Let $G = \{1234\}$ and let $Q$ be the set $q_{12}q_{21} =
q_{34}q_{43} = 1$. The associated space of constants is generated
by
$
  e_1e_2 - q_{21}e_2e_1$ and $e_3e_4 - q_{43}e_4e_3 .
$
One finds that
$$
\det S_G(Q) =
(1-\sigma_{13})^6(1-\sigma_{16})^6(1-\sigma_{23})^6(1-\sigma_{24})^6
(1-\sigma_{1234}).
$$
There is a 1-dimensional subspace of ${\cal B}_G(Q)$ of primitive
constants, on the primitive boundary of $V^Q$ on which
$\sigma_{1234} = 1$.
\b

\no{\bf 2.3.5.   Classification by multilinear type in
Case $N = 3$.}

The constraints of order 2 are, up to permutations of the indices,
$$
Q^1: \sigma_{23} = 1,~Q^{12}: \sigma_{13} = \sigma_{23} = 1,~
Q^{123}: \sigma_{23} = \sigma_{12} = \sigma_{13} = 1,
$$
and the
empty set (no constraint). The classification by multilinear
types of total degree 2 yields 4 types. The discussion of Example
2.3.3 shows  that there is a distinguished boundary only if
$Q$ is the empty set. The complete classification by multilinear
type thus yields 5 distinct types (up to a permutation of the
generators).
  \b

\no{\bf 2.3.5.   Classification by multilinear type in
Case $N = 4$.}

At total order 2 there are  11 possibilities (always
up to a permutation of the generators). We list the set of
parameters that are fixed at unity in each case.
\b
1. None,

2. $\sigma_{12}$,

3. $\sigma_{12}, \sigma_{34}$,

4. $\sigma_{12}, \sigma_{13}$,

\vskip-1.2cm\hskip1in 5. $\sigma_{12}, \sigma_{13}, \sigma_{14}$,

\hskip1in 6. $\sigma_{12}, \sigma_{13}, \sigma_{23}$,

\hskip1in 7. $\sigma_{12}, \sigma_{23}, \sigma_{34}$,

\vskip-1.65cm\hskip2.2in 8. $\sigma_{12} , \sigma_{23},
\sigma_{34},
\sigma_{14}$,

\hskip2.2in 9.  $\sigma_{12}, \sigma_{23},
\sigma_{34}, \sigma_{13}$,

\hskip2.2in 10.  $\sigma_{12}, \sigma_{23},
\sigma_{34}, \sigma_{13}, \sigma_{14}$,

\hskip2.2in 11. all  $\{\sigma_ {ij}\}_{i<j}$.
\b
\no This give rise to 11 varieties
$V^{Q_1}, ... , V^{Q_{11}}$. A distinguished boundary at degree 3
appears in the first five cases only. Up to total order 3 there are
16 possibilities: the 11 cases listed and in addition the
following,
\b
12. $\sigma_{123}$,

13. $\sigma_{123},\sigma_{124}$,

14. $\sigma_{123},\sigma_{124}, \sigma_{134}$,

15. all $\{\sigma_{ijk}\}_{i<j<k}$,

\vskip-1.2 cm\hskip1.5in 16. $\sigma_{12},\sigma_{134}$,

\hskip1.5in 17. $\sigma_{12}, \sigma_{14},
\sigma_{234}$,

\hskip1.5in 18. $\sigma_{12}, \sigma_{13},
\sigma_{14}, \sigma_{234}$,

\b
Finally, we examine each of $V^{Q_1}, ... ,V^{Q_{16}}$ and find
that there is a primitive boundary in cases 1,2,3 and 12 only. the
complete list of multilinear types for $N = 4$ is given by the 16
possibilities already listed, plus the following,
\b
19.  $\sigma_{1234}$,

20. $\sigma_{12}, \sigma_{1234}$,

\vskip-8mm \hskip1.5in 21. $\sigma_{12},\sigma_{34},
\sigma_{1234}$.

\hskip1.5in 22. $\sigma_{123}, \sigma_{1234}$.
\bb
\topskip-.5in

\no{\bf 2.3.6. The general case.}

Our problem  can be solved by calculating the
determinant of
$S_G(Q)$ for all degrees $G$ and for any set $Q$ of
constraints $\sigma_{i_1\cdots i_k} = 1, ~k <n$. Eq.
(2.2) gives the answer in the simplest case, when there
are no constants of degree lower than
$G$. Our first result, Theorem 4.2, gives the condition for $V^Q$
to have a primitive boundary and for the existence of a primitive
constant of degree $\{1,2,...,n\}$ in the case that
all primitive constants of lower degree have
total degree 2. The most far reaching result obtained is
Theorem 6.5, which applies whenever there is at
least one $\sigma_{ij}$ that is not subject to any constraint.
\bb
\no{\steptwo 3.  Varchenko's method. }
\b
\ce{\bf 3.1  Arrangements of hyperplanes.}

Following
Varchenko [V], we consider an arrangement of hyperplanes
$H_1,\cdots , H_k$ in $\Rit^n$. An edge is a non-empty
intersection of hyperplanes and a domain is a connected
part of the complement of the set of hyperplanes.
Complex weights $a_1,a_2,\cdots $ are attached to the
hyperplanes, the weight $a(L)$ of an edge $L$ is the
product of the weights of the hyperplanes that contain
$L$. A bilinear form is defined by
$$
B(D,D') = \prod a_i,
$$
where $D,D'$ are any two domains and the product runs
over the hyperplanes that separate them. Varchenko gives
the following formula [V],
$$
\det B = \prod_L\bigl(1-a(L)^2\bigr)^{n(L)p(L)}.\eqno(3.1)
$$
The product runs over all edges and  $n(L), p(L)$ are
certain natural numbers or zero.

\b
\ce{\bf 3.2. Interpretation.}

Let the parameters
$\{q_{ij}\}$ be in general position. For
a special choice of hyperplanes and weights, $B$ is
identified with the matrix $S_G,~ G = \{1,\cdots,n\}$.
Namely, let the hyperplanes be
$$
H_{ij} = \{x_i = x_j\},\quad i<j,\quad i,j = 1,\cdots,n.
$$
The domains are then
$$
D_{\underline i} = \{x_{i_1}< x_{i_2} < \cdots
< x_{i_n}\},
$$
in natural correspondence with the monomials
$e_{i_1}\cdots e_{i_n}$ of degree
$G = \{i_1,\cdots,i_n\}$. Choose $q_{ij}$ symmetric (this
is harmless since the zeros of $\det S_G$ depend only on
the products
$q_{ij}q_{ji} = \sigma_{ij})$ and set $q_{ij} =
a_{ij} $,   the weight of the hyperplane
$H_{ij}$. Then
$$
B(D_{\underline i},D_{\underline j}) = S_{\underline
i\underline j}.
$$

For this case, Varchenko's   formula (3.1) for $\det B$
   coincides with the Eq.(2.4).
The contributing edges are all those of the form $L =
\{x_{i_1} = x_{i_2} = \cdots =x_{i_k}\}$. In this\break
interpretation the integer
$n(L)$ is the number of domains in the arrangement in
$L$ for\break
\no  which
  the hyperplanes    are the
intersections of
$L$ with those of the original planes that do not contain
$L$. This number is the same as the number $u(n,k) =
(n+1-k)!$ in Eq.(2.2); that is, the dimension of the
intersection between
$\bbb_{G}$ and the ideal generated  by any  primitive
constant of degree $\{1,...,k\}$.
\b
\ce{\bf  3.3. The number $p(L)$.}

The calculation of the
number
$p(L)$, in the exponent in Eq.(3.1), is more subtle and at
the center of interest. Any edge of the configuration
under consideration is a hyperplane of the form (up to a
renaming of the coordinates)
$$
\{x_1 = \cdots =x_k,~ x_{k+1}=\cdots =
x_{k+l },~\cdots~,~ \cdots \}.
$$
   The number $p(L)$ is defined as follows.
Let $N$ be the normal to $L$; in our case it is
$$
N = \{ \xi_1,\cdots,\xi_k;~ \eta_1,\cdots
,\eta_l;~\cdots~;~\cdots\},\quad \sum \xi_i = \sum \eta_j
=
\cdots = 0.
$$
There is the arrangement $\{H\cup N; L\subset H\}$~of
hyperplanes in $N$. The planes are
$\xi_i = \xi_j, ~ \eta_i = \eta_j,~\cdots ~ $. Consider
the projectivization of this arrangement. Fix any
one of the hyperplanes, $H$,~say.~Then the number
$p(L)$ is the number of projective domains the closures
of which do not intersect $H$.
\b
\no{\bf Proposition.} {\it The number $p(L)$ is zero unless
$L = \{x_1 = \cdots =x_k\}$ for some $k = 2,\cdots,n$,
up to a permutation of the index set.}
\b
\no{\bf Proof.}
Let $H = \{\xi_1 =
\xi_k\}$. If the set of coordinates of $N$ includes one or more
sets beyond the initial set
$\xi_1,\cdots,\xi_k$, then all the projective domains
include   some lines on which $\xi_1 = \cdots = \xi_k
= 0$, and then $p(L) = 0$.   The proposition is proved.
\b
\no{\bf Corrollary.} {\it The determinant of $S_G, ~G =
\{ 1,\cdots n\}$, is a product of factors of the form}
$1-\sigma_{i_1...i_k},~ k = 2,\cdots ,n$.
  \b
We are reduced to the case when, up to a renaming of the
coordinates,
$$
L = \{x_1 = \cdots =
x_k\},   \quad  N =
\{\xi_1,\cdots \xi_k,0,\cdots ,0\},~~ ~\sum_i^k\xi_i = 0.
\eqno(3.2)
$$
\no Remember that the parameters are in general position, no
constraints.

Consider the closure  $\xi_1 \leq \cdots
\leq \xi_k$ of the domain $\xi_1 < \cdots < \xi_k,~ ~\sum
\xi_i = 0$. It touches any hyperplane
$\{\xi_i = \xi_j\}$ at points where $ \xi_i = \cdots
= \xi_j = 0$. What saves us from the
conclusion that $p(L)$ is always zero is the fact that
this domain fails to touch the hyperplane $\xi_1 =
\xi_k$. This is because the point at the origin of $N$
does not have a projective image. We conclude that for
such edges, $p(L) = (k-2)!$, the number of domains of
the type $\xi_1 \leq\cdots \leq \xi_k$. This number is
the same as the number $v_k = (k-2)!$ in (2.3);
that is, the dimension of the space of primitive constants of
  degree $G_k$, so
Varchenko's formula (3.1) reduces to Eq.(2.4) in this
case.

\bb

\no {\steptwo 4. Constants of lower degree, each of
total  order 2.}
\b
\ce{\bf 4.1. A special case.}

We shall determine under
what conditions there are primitive constants of degree
$G =
\{1,\cdots,n\}$ in the case that there is any number of
primitive constants of lower degree, but all of them of
total order 2. The parameters are thus in general
position, except that they satisfy a set of constraints,
$$
\sigma_{ij} = 1,\quad \{i,j\} \in P,
$$
where $P$ is a fixed subset of the set of
pairs $\{i,j\}, ~i\neq j ,~ i,j = 1,\cdots,n$.
  The primitive constants of lower degree are $e_ie_j -
q_{ji}e_je_i,~
\{i,j\} \in P$.

In the idiom of arrangements of hyperplanes, the
constraint $\sigma_{ij} = 1$ means that the weight
$a_{ij}$ of the plane $x_i = x_j$ is equal to unity.

Let $B(Q)$ be the matrix $B$ for the arrangement obtained by
removing hyperplanes with weight 1. The identification $S_G(Q) =
B(Q)$ holds in this case as well. The arrangement is in
$\Rrm^n$, with hyperplanes
$x_i = x_j, ~ \{i,j\} \notin P$.

Varchenko's formula (3.1) applies and the
only edges for which the number $p(L)$ is different from
zero   are the ones of the form (3.2). The determinant
still has the form (2.4), but some of the exponents are
diminished. Of importance for the classification problem
is  the question whether the factor
$1-\sigma_{1\cdots n}$ appears with non-zero exponent: if $L_0
:=
\{x_1 = \cdots = x_n\}$,
when is $p(L_0)$ different from zero?
\b
  The space $N$ normal to $L_0$ is
$\{\xi_1, \cdots ,\xi_n\},~ \sum \xi_i = 0$. The domains
are defined by inequalities,
$$
\xi_i <\xi_j,\quad \{i,j\} \notin P.\eqno(4.1)
$$
A domain the closure of which does not intersect a given
hyperplane, $\xi_1 = \xi_n$, say, must bracket all the
other variables, $\xi_2,\cdots ,\xi_{n-1}$, between
$\xi_1$ and $\xi_n$.
\b
\ce{\bf 4.2. Theorem.}

{\it Let $G = \{1,\cdots,n\},~ n > 2$, $Q$ the set of constraints
$\sigma_{ij} = 1, \{i,j\} \in P$.
We may suppose that $\sigma_{1n} \neq 1$; that is, that
$\{1,n\} \notin P$.   The following condition
is necessary and sufficient for $V^Q$ to have a primitive boundary.
For any
$i, ~1<i<n$, there is a sequence
$1,\cdots,i,\cdots ,n$, a subsequence of a permutation
of
$ 1,\cdots ,n$, such that no pair of neighbours in it
belongs to $P$.}

  \b

\no{\bf Examples.}  See the list of constraints in 2.3.5. In the
case that $P = \{(12), (34)\}$ the projective domain $1<3<2<4$ does
not touch the plane $x_1=x_4$ since $P$ does not contain the
pairs (1,3), (3,2) or (2,4). There is only one such domain, so the
determinant contains  the factor $1-\sigma_{1234}$ with exponent
1.
   There is at least one
(actually exactly one) primitive constant of degree $G =
\{1,2,3,4\}$ when $Q$ is the set
$\sigma_{12} = \sigma_{34} = 1$. But if the
constraint is
  $\sigma_{12} =
\sigma_{13}=1$, then there is no primitive constant of this degree.
  See [Kh3].

\bb
\no{\steptwo 5. Algebraic proof of Theorem 4.2.}
\b
It was seen that the case of constraints of a very
special type lies within the range of the theory of
arrangements of hyperplanes. But the direct application
of this theory to more general situations does not appear to be
straightforward. For that reason it will be useful to reformulate
the proof of Theorem 4.2 in purely algebraic terms.

By stipulation, the relations of $\bbb_q(Q)$ are generated by
$e_ie_j = q_{ji}e_je_i,~\{i,j\}\in P$, and the constraints are
$Q: \sigma_{ij} = 1, ~ \{i,j\} \in P$. We may suppose that
$\{1,n\} \notin P$.
\b

\ce{\bf 5.1. Basis.}

All bases used for $\bbb_G(Q)$
and its subspaces will be monomial. A basic monomial will be
called a word. If
$e_{\underline i}$ is a word then so is
$e_{\underline i'}$ (the same word read backwards), unless the two
are proportional to another. Since $\sigma_{1n} \neq 1$, this
cannot happen in the context.
\b
  \no{\bf   Example.} If $\sigma_{12} = 1$, then a basis for
$\bbb_{\{1,2,3\}}(Q)$ is
$$
e_1e_3e_2,\quad e_2e_1e_3,\quad e_2e_3e_1,\quad
e_3e_1e_2.
$$
We shall say that a word in $\bbb_G(Q)$  is
`positive' if
$e_1$ precedes $e_n$, and proceed to choose the positive words of
a basis.

\b
\ce{\bf 5.2.   Factors and classes.}

A positive word
in $\bbb_G(Q)$ has the form $xe_1ye_nz$. The degree
  of $y$ defines a filtration of $\bbb_G(Q)$. Choose
a monomial basis that respects this filtration. An
element of the basis, of the form $xe_1ye_nz$, will be
said to have the `factor' $y$ and to be of `class' $g$ =
the degree of $y$. Remember that the degree of $e_{\underline i}$
is the unordered set $\{i_1,...,i_k\} $ of indices.

\b
\no{\bf  Lemma.} {\it The number of words with
factor $y$ depends only on the degree  of $y$.}
\b
\no{\bf Proof.} The degree of $y$ selects a subset of the
generators $e_2,...,e_{n-1}$ and reduces the construction of the
basis of ${\cal B}_G(Q)$ to that of a basis for the subspace  of
polynomials in $e_0:= e_1ye_n$ and the supplementary set of
generators. The relations are all of the type $e_ie_j = ke_je_i,~ k
\in \Crm$ and are independent of the order of factors in the
monomial $y$.

\b
\ce{\bf 5.3. First sum rule.}

Let $u_n(g)$ be the number
of positive words for $\bbb_G(Q)$ that contain some fixed
factor $y$ of degree   $g$, and let $v(g)$ be the number of
(linearly independent) factors of degree $g$. Then
$$
\sum_g u_n({  g})v({  g})  = \1/2 \dim
\bbb_G(Q).
$$

  \b
\no {\bf   Example.}  $P =
\{\{1,2\}\}$.
$$
\matrix{Type & u &v  &{\rm Basis,~positive ~part~},
i\neq j = 2,3\cr\cr (\cdot) & 6 &1&
e_ie_j(e_1e_4),~~e_i(e_1e_4)e_j,~~(e_1e_4)e_ie_j\cr
(3)& &0\cr (4) &2   & 1&e_2(e_1e_3e_4),~
(e_1e_3e_4)e_2\cr (34)&&0&\cr
(32)&1&1 &(e_1e_3e_2e_4)
\cr}
$$
\b
\ce{\bf 5.4. Second sum rule.}

\no{\bf  Lemma.} {\it Let $\tilde v(g)$ be the exponent
of $(1-\sigma_{1  {\underline i} \,n})$ in $\det S_{\hat
  g}(Q)$,
$\hat g = \{1, g,n\}$,
$g$ the degree of $e_{\underline i}$; then the
exponent of the same factor in $\det S_G(Q)$ is $ u_n(g)\tilde
v(g)$.}
\b
\no{\bf Proof.} Fix $\hat g = \{1,g,n\} =\{1,
i_2,\cdots,i_k,n\}$. In the matrix $S_G(Q)$ replace by zero all
$\sigma_{ij},\, \{i,j\}
\notin P$, that do not appear in $\sigma_{1\underline
i\,n}$. Then $ S_G(Q)_{\underline j\underline k}$ vanishes
unless $e_{\underline j}$ and $e_{\underline k}$ are
equal up to a reordering of the generators
$e_1,e_{i_2},...,e_{i_k},e_n$ only; that is, unless $e_{\underline
j} = xe_1ye_nz,
\, ~e_{\underline k} = xe_1y'e_nz$ with $y$ and $y'$ of  the
same degree
$g$. The matrix takes the block form and $\det S_G(Q)$ reduces to a
power of $\det S_{\hat g}(Q)$. The exponent is the number of
blocks
and is equal to the number
$u_n(g)$ of words that contain some fixed $y$
of grade
$g$.

\b

Every $\sigma_{1 \underline
i\,n}$ is  linear in $\sigma_{1n}$, and the terms of highest power
of $\sigma_{1n}$ in $\det S_G(Q)$
$$
\prod \partial_{\underline i'} e_{\underline i}
\propto (\sigma_{1n})^\kappa,\quad \kappa = \1/2 \dim
\bbb_G(Q),
$$
where the product runs over all words. Hence
$$
\sum_g u_n(g)\tilde v(g) = \1/2 \dim \bbb_G(Q).
$$
Clearly, $\tilde v(g) = v(g)$ when $n = 1$; therefore by
induction in  $n$, it follows from the two sum rules that
$\tilde v(g) = v(g)$ and Theorem 4.2 is proved (again).
\b\b
\no{\steptwo 6. General constraints.}

\ce{\bf 6.1. A proviso.}

So far we have allowed primitive relations of
order 2 only, associated with constraints $\sigma_{ij} =
1,\,\,\{i,j\} \in P$, where $P$ is a collection of
pairs. The key to both proofs of Theorem 3.2 was the
counting of powers of $\sigma_{1n}$ in $\det S_G(Q)$, and
for this reason it was essential that $\{1,n\} \notin
P$. This is not a real limitation, for another pair will
do just as well, as long as there is at least one that
is not in $P$; that is, except in the case that $\det
S_G(Q) = 1$. Now let us consider the more general
situation, when there is a family of constraints,
$$
\sigma_{\underline i} = 1,\quad \underline i \in
P,\eqno(6.1)
$$
where $P$ is any collection of proper subsets of
$\{1,\cdots, n\}$. To apply the method of counting
powers of $\sigma_{1n}$, we need for this parameter to
be unconstrained, and this amounts to the limitation
that, for any $g$,
$$
\{1,g,n\} \notin P.\eqno(6.2)
$$

   Proceeding as in
Section 5, we encounter no difficulties in choosing a monomial
basis based on the concepts of `factors' and `classes'. But the
Lemmas in 5.2 and 5.4  need  to be re-examined. The first one is:
\b
\no{\bf  6.2. Lemma.} {\it The number of words with
factor $y$ depends only on the degree  of $y$.}
\b
Let us
consider the process of choosing the basis in somewhat more
detail. Begin with elements of the form $e_1ye_n$, $y$ a
permutation of $e_2
\cdots e_{n-1}$. If there are no relations involving either $e_1$
or
$e_n$, then any set of independent `$y$'s will do.
It is enough to consider primitive constants involving
$e_1$, say. (By stipulation, there are no primitive constants
involving both $e_1$ and $e_n$.) Any constant of this type is a
polynomial
$$
e_1A + \sum_{j = 2}^{n-1} e_je_1A^j + \cdots + Be_1.
$$
If there is only one such constant, with $A \neq 0$, then
the filtration replaces $y$ by $y$ modulo the right ideal
$\{Az\}$. This affects the number $v(g)$ of factors of
this type, but the relation has no further effect on the
construction of the basis, and the lemma still stands.
If there is another constant, of the same type, then the
argument still applies. But   it may happen
that there are two constants with the same A, and then
there is a constant of the type
$$
\sum_je_je_1A^j + \sum_{j,k} e_je_ke_1A^{jk} + \cdots +
Be_1.
$$
In this case, suppose $A^2 \neq 0$, then differentiation with
   $\partial_2$ gives
$$
e_1A^2 + \partial_2\bigl(\sum_{j>3}e_je_1A^j +
\sum_{j,k}e_je_ke_1A^{jk} +
\cdots + Be_1\bigr) = 0.
$$
Now this relation implies that the factor $y$ is defined
modulo the ideal $\{A^2z\}$. Continuing in this manner
we conclude that any relation that involves $e_1$ leads
to a reduction in the number of factors, but it does not
affect the number of basis elements containing a given
factor. In fact, once the factors have been determined,
then the enumeration of basis vectors is independent
of the factor and depends only on the set of generators
in it; that is, on its degree.   The Lemma is proved.

This
implies that the first sum rule remains valid.

\b
\no{\bf Example.} $P = \{\{1,2,3\}\}$, thus
$\sigma_{123} = 1, n = 4$.
$$
\matrix{ {\rm Factors} & u_4 & v & {\rm words}\cr (\cdot )
&6&1&e_ie_j(e_1e_4), ~ e_i(e_1e_4)e_j,~ (e_1e_4)e_ie_j\cr
(e_2)&2&1&e_3(e_1e_3e_4),~ (e_1e_2e_4)e_3\cr
(e_3)&2&1&e_2(e_1e_3e_4), ~(e_1e_3e_4)e_2\cr
(e_2e_3)&1&1&(e_1e_2e_3e_4)
\cr}
$$
There is only one constraint and its only effect is to exclude
the monomial $(e_1e_2e_3e_4)$ with factor $e_3e_2$ from the
basis.
\b

The second sum rule also remains in force in the more
general case, but our proof of the Lemma in 5.4 does not, since
we cannot replace by zero parameters that are
constrained. It must be replaced by the following two lemmas.

Fix a natural number $k,~ 1<k <n$ and a subset $\{i_2,...,
i_k\}\subset G = \{1,...,n\}$. Fix a set
$Q$ of constraints so as to leave the parameter
$\sigma_{1n}$ free, and let
${\cal B}(Q)$ be the associated quotient algebra.   Write
   $\{i_2,\cdots ,i_k\} = g$.
\b
\no{\bf 6.3. Lemma }
{\it Fix all parameters except $\sigma_{1n}$. The exponent of
$(1-\sigma_{1i_2...i_kn})$  in $\det S_G(Q)$ is equal to the
dimension of $\ker S_G(Q)$ at the value of $\sigma_{1n}$
that makes $\sigma_{1i_2...i_kn} = 1$.}

  \b
\no{\bf Proof.}  Introduce a monomial basis as
above, consisting of `positive' words in which $e_1$
precedes $e_n$, and the same set taken in reverse order.
Set $q_{1n} = q_{n1} = q$. Then the matrix element $S_{\underline
i\underline j}(Q)$ of $S_G(Q)$  is independent of $q$ if
$\underline  i,\underline j$ are both positive or both
negative,   linear in $q$ otherwise, and $S_G(Q)$ takes
the form
$$
S_G(Q)  = \pmatrix{A&qB^t\cr qB & C\cr},
$$
with  $A,C$
symmetric (take $q_{ij} = q_{ji}$) and invertible(there are no
constants in the subalgebra of ${\cal B}_G$ generated by
$e_2,...,e_{n-1}$). Interpreting
$S_G$ as a form, we transform $A\oplus C$
to a unit matrix without affecting the kernel of
$S_G(Q)$, converting this form to  $I + qD$ with $I_{ij} =
\delta_{ij}$ and
  $D$ symmetric. Finally, interpreting $I + qD$ as a symmetric (and
hence diagonalizable) matrix  one obtains the result.

\b
\no{\bf 6.4. Lemma.}
{\it If
$(1-\sigma_{1i_2...i_kn})$ appears with exponent
$\tilde v(g)$ in $\det S_{\hat g}(Q)$, $\hat g = \{1, g,n\}$,
then it appears with exponent $u_n(g)\tilde v(g)$ in
$\det S_G(G)$.}
\b
\no{\bf Proof.}  If $xe_1ye_nz$ is a word of
class $g = \{i_2,..., i_k\}$, and $C \in \bbb_{\{1,g,n\}}(Q)$ is
a constant, then $xCz$ is in the ideal generated by $C$, and
this correspondence extends to a bijection (for fixed $y$ and
$C$). The dimension of the ideal generated in ${\cal B}_G(Q)$ by
the constants in
${\cal B}_g(Q)$ (which is the same as the dimension of $\ker
S_G(Q)$ at $\sigma_{1i_2...i_kn} = 1$,
and by Lemma 6.1 equal to the exponent of $1-\sigma_{1i_2...i_kn
}$ in $\det S_G(Q)$) is thus
$u_n(g)$ times the dimension of the space of constants in ${\cal
B}_g$. The lemma is proved.
\b
This gives the second sum rule, and by induction,
   $v(k) = \tilde v(k)$. We have thus proved:
\b
\no{\bf 6.5. Theorem.} {\it Assume parameters as before, with Eq.s
(6.1) and (6.2). Then the factor $(1-\sigma_{1i_1... i_k,n})$
appears in $\det S_G(Q)$ with an exponent that is equal to
the number of words of class
   $g = \{2,\cdots ,k\}$ in ${\cal B}_G(Q)$  and
   the number of linearly independent
constants that appear in
${\cal B}_G(Q)$ when $\sigma_{1\cdots n}$ tends to 1 is equal to
the number of words of class $\{2, \cdots , n\}$.

Cosequently, a necessary and sufficient condition for the
appearance of at least one primitive constant in ${\cal B}_G(Q)$
when $\sigma_{1\cdots n}$ tends to 1, and for the existence
of primitive boundary of $V^Q$, constant is that there is
$e_1ye_n \in \bbb_G(Q)$ that cannot be expressed in terms
of monomials $xe_1y'e_nz$ with $xz$ of non-zero degree.}

\b
It remains to understand the case when all the
parameters are constrained. This can happen with as few
as 3 independent constraint, for example $\sigma_{1\cdots
n-1} = \sigma_{2\cdots n} = 1$ and $\sigma_{1n} = 1$.
Present methods fail because it makes no sense to count
the powers of any one of the parameters.

\bb

\ce{\steptwo Acknowledgement.}
\b
I thank Vladislav Kharchenko, Geeorge Pinczon and Daniel
Sternheimer for helpful comments.
\ve
\ce {\steptwo References.}
\ss
\no [D] V.G. Drinfel'd, ``Quantum Groups", Proc. ICM
(Berkeley 1986). Amer. Math. Soc. \no

\quad 1987.
\smallskip
\no [F1] C. Fronsdal, Generalizations and deformations of
quantum groups, Publ. Res. Inst.

\quad Math. Sci. {\bf 33},
91-149 (1997). (q-alg/9606020)
\smallskip

\no [F2] C. Fronsdal, On the classification of q-algebras,
Lett. Math. Phys. {\bf 222}, 708-746

\quad (1999).
\smallskip

\no [FG] C. Fronsdal and A. Galindo, The ideals of free
differential algebras, J. Algebra {\bf 222}

\quad (1999) 708-746.
\smallskip

\no [FCG] Flores de Chela D. and Greene J.A., Quantum Symmetric
Algebras, Algebras and

\quad Rep.s Th.
{\bf 4} 55-76 (2001).
\smallskip

\no [Ka]~~~ M. Kashiwara, On crystal basis of the Q-analogue of
universal enveloping algebras,
\smallskip
\quad Duke Math. J. {\bf 63},
465-516 (1991).
\smallskip

\no [K]~~~~ V. Kac, {\it Infinite dimensional Lie algebras,
Cambridge University Press 1990.}
\smallskip

\no [Kh1]~ V. Kharchenko, Differential calculi and skew
primitively generated Hopf algebras,

\quad preprint.
\smallskip

\no [Kh2]~ V. Kharchenko, An Existence Condition for
Multilinear Quantum Operations, J.

\quad Algebra {\bf217},
188-228.
\smallskip
\no [Kh3]~ V. Kharchenko,  Skew primitive elements in Hopf
algebras and related identities,

\quad Journal of
Algebra {\bf 238}(2001) 534--559, section 7;
    An algebra of skew primitive

\quad elements,
       Algebra and Logic {\bf 37}
N2 (1998) 101-126, theorems 8.1 and
   8.4
\smallskip

\no[R]~~~~~M. Rosso M., Quantum groups and quantum shuffles,
Invent. Math.{\bf 133} (1998)

\quad 399-416.
\smallskip
\no [V]~~~~ A. Varchenko, Bilinear Form of Real Configuration
of Hyperplanes, Adv. Math. {\bf 97},\break

\vskip-4mm\quad 110-144
(1993).
Multidimensional Hypergeometric Functions,
Representation\

\quad  Theory of Lie algebras and
  Quantum
Groups, Adv. Ser. Math. Phys. Vol. {\bf 21}.

\quad World
Scientific, Singapore, 1995.
  Adv.Math. {\bf 97}, 110-144
(1993).
\smallskip

\no[Va] ~~~V.S. Varadajan,
\smallskip

\no [W]~~~~L. Woronowicz, Differential calculus on compact matrix
pseudogroups (quantum

\quad groups), Comm. Math.
Phys. {\bf
122}, 125-170 (1989).
\smallskip

\no [Y]~ ~~ H. Yamane, A Serre Type Theorem for Affine Lie
Algebras and their Quantized

\quad Enveloping Algebras,
Proc. Japan Acad. {\bf 70}, ser. A  31-36 (1994).
\end

\end